\documentclass[11pt]{article}

\newtheorem{theorem}{Theorem}
\newtheorem{lemma}[theorem]{Lemma}
\newtheorem{corollary}[theorem]{Corollary}
\newtheorem{proposition}[theorem]{Proposition}
\newenvironment{proof}{{\it Proof:}\hspace {0.5em}}
{\\* \null\hfill  \par}

\newenvironment{exafont}{\begin{bf}}{\end{bf}}
\newenvironment{algorithm}{\vspace{0.3cm}\par\noindent\refstepcounter{theorem}
\begin{exafont}Algorithm \thetheorem\end{exafont}\hspace{\labelsep}}{\vspace{0.3cm}\par}
\newenvironment{problem}{\vspace{0.3cm}\par\noindent\refstepcounter{theorem}
\begin{exafont}Problem \thetheorem\end{exafont}\hspace{\labelsep}}{\vspace{0.3cm}\par}

\newfont{\Bb}{msbm10 scaled\magstep 1}

\newcommand{\conv}{\mbox{conv}}

\newcommand{\R}{\mbox{\Bb R}}
\newcommand{\Z}{\mbox{\Bb Z}}

\title{\Large  Decomposition of Polytopes and Polynomials
\footnote{
The first author was supported in part by NSF under Grant \#DMS9970637 and 
NSA under Grant \#MDA904-00-1-0048.
The second author gratefully
acknowledges the support of the Marr Educational Trust and
Wolfson College, Oxford, and thanks Dominic Welsh for his help
and encouragement.}
}

\author{Shuhong Gao\footnote{
Department of Mathematical Sciences,
Clemson University, Clemson, SC 29634-0975, USA.
E-mail: sgao@math.clemson.edu.}
and Alan G.B. Lauder\footnote
{Mathematical Institute, Oxford University, Oxford OX1 3LB, U.K.
E-mail: lauder@maths.ox.ac.uk.
}}
\date{November 3, 2000}

\begin{document}

\baselineskip = 1.5\baselineskip 

\maketitle

\begin{abstract}
Motivated by a connection with the factorization of multivariate
polynomials, we study integral convex polytopes and their integral
decompositions in the sense of the Minkowski sum. 
We first show that deciding decomposability of integral polygons is
NP-complete then present a pseudo-polynomial time algorithm for
decomposing polygons.  For higher dimensional polytopes, we give 
a heuristic algorithm which is based upon projections and uses randomization.  
Applications of our algorithms include absolute irreducibility testing
and factorization of polynomials via their Newton polytopes.
\end{abstract}

\section{Introduction}\label{sec-intro}

It is well-known that the theory of convex polytopes has many
applications across mathematics and computer science 
\cite{Ewa96,GO97,GW93,Stu96}.
One such application is to polynomial factorization, and
motivated by this connection we discuss decomposition algorithms
for polytopes.
Given a multivariate polynomial one may
associate with it, in a way we shall fully explain in
Section \ref{sec-prel}, an integral polytope called its Newton polytope.
It was observed by Ostrowski in 1921 that
if the polynomial factors
then its Newton polytope decomposes, in the sense of the
Minkowski sum, into the Newton polytopes of the factors.
The ramifications of this simple observation are two-fold. Firstly, 
criteria which ensure polytope
indecomposability can be used to construct families of irreducible,
indeed absolutely irreducible, polynomials. Secondly, algorithms which
test whether a polytope is decomposable and construct
decompositions may be useful in factoring polynomials. Of course, 
such criteria and algorithms are also of independent interest
and may have other applications. Indecomposability conditions were
explored by the first author in \cite{Gao99} and will be discussed further
in Section \ref{sec-indecomp}.
Our main focus will be, however, on the second application, that is
on algorithms for decomposing polytopes.

We first show that the problem of
testing whether a polytope is indecomposable is NP-complete even in 
dimension two, so there does not exist, unless NP $=$ P,
a genuinely efficient algorithm for decomposing polytopes.
However, we present a ``pseudo-polynomial'' time algorithm 
(see \cite{GJ79}) for testing 
indecomposability in dimension two and a modified version which also
allows one to count the number of decompositions and find summands.
We also discuss a heuristic algorithm which uses randomization for
testing higher dimensional polytopes for indecomposability. In Section
\ref{sec-appl}, we describe applications of our algorithms to
polynomials with respect to their irreducibility and factorization.
In particular, we touch upon 
an open
problem in polynomial factorization which we now describe.
In his survey paper on polynomial factorization \cite{Kal92},
Kaltofen concludes with several open problems
one of which, due to B. Sturmfels, is stated as follows:
``From the support vectors $(e_{j 1}, \ldots, e_{j n})$
of a sparse polynomial 
$ \sum_{j=1}^{t} a_{e_{j 1},\ldots, e_{j n}} X_1^{e_{j 1}} \cdots 
   X_n^{e_{j n}}$,
compute by geometric considerations the support vectors of
all possible factorizations''.
This problem can be attacked
by our polytope method, although it must be noted that we are unable
to give a complete solution.
The basic idea runs as follows:
Given a bivariate polynomial, we can compute its Newton polytope
and then find all the integral summands of this polytope.
The summands correspond to the Newton polytopes
of all the possible factors of the polynomial.
The integral points in a summand
give the support vectors of the factor corresponding to the summand.

The remainder of the paper is organized in the following way.
Section \ref{sec-prel} contains the necessary background material
on the theory of convex polytopes and in Section \ref{sec-indecomp} 
we discuss some preliminary results on
polytope indecomposability which shall be useful to us but are also
of independent interest. Section \ref{sec-decomp} is devoted to
algorithms and is further divided into two parts: In Section \ref{sec-polygon}
we present algorithms
for both testing polygons for decomposability and counting and
constructing decompositions of polygons. Section \ref{sec-polytope} contains a
heuristic randomized algorithm for higher dimensional polytopes based
upon projections down to dimension two. Finally, in Section \ref{sec-appl} we
discuss applications of these algorithms to absolute irreducibility
testing and polynomial factorization.

\section{Polynomials and  Newton polytopes}\label{sec-prel}

\subsection{Background geometry and algebra}

Before describing the connection between polynomials and
polytopes, we recall some terminology and results from
the theory of convex polytopes (\cite{Gru67}). Let $\R$ denote the field of
real numbers and $\R^n$ the Euclidean $n$-space. A {\em convex set} in
$\R^n$ is a set such that the points on the line segment
joining any two points of the set lie in the set; the {\em convex
hull} of a set of points is the smallest convex set which contains
them; and  the convex hull of a finite set of points is called a {\em
convex polytope}. A point of a polytope is called a {\em vertex} (or
{\em extreme point})  if it does not belong to the interior of 
any line segment contained in the polytope. A polytope is always
the convex hull of its vertices.  A hyperplane {\em cuts} a polytope if both of the open half spaces
determined by it contain points of the polytope. A hyperplane which
does not cut a polytope, but has a non-empty intersection with it is
called a {\em supporting hyperplane}. The intersection of a supporting
hyperplane and a polytope is a {\em (proper) face}, and the union of
all (proper) faces is the {\em boundary}. One may equivalently define
a vertex to be a $0$-dimensional face, and $1$-dimensional faces are
known as {\em edges}.

For two subsets $A$ and $B$ in $\R^n$, define their
{\em Minkowski sum} to be $A+B = \{a+b\,|\,a\in A,b\in B\}$. We call 
$A$ and $B$ the {\em summands} of $A+B$. It is easy to show that the
Minkowski sum of two convex polytopes is a convex polytope.

Let $f\in K[X_{1},\ldots,X_{n}]$ be a nonconstant polynomial
where $K$ is an arbitrary field.
We call  $f$ {\em absolutely irreducible} over $K$
if it has no non-trivial factors over the algebraic closure of $K$.
Suppose
$$ f = \sum a_{i_{1}\ldots i_{n}} 
  X_{1}^{i_{1}}\cdots X_{n}^{i_{n}}.$$ 
For each term with $a_{i_{1}\ldots i_{n}} \ne 0$, the 
corresponding exponent vector $(i_{1},\ldots,i_{n})$, viewed
in $\R^n$, is called a support vector of $f$. Define $Supp(f)$
to be the set of all support vectors of $f$, i.e.,
\[ Supp(f) = \{(i_{1},\ldots,i_{n})\,|\,a_{i_{1}\ldots i_{n}} \ne 0\}.\]
Note that $Supp(f)$ is empty if $f=0$.
The {\em total degree} of $f$, where $f \ne 0$, is
the maximum value of $\sum_{1 \leq j \leq n} i_{j}$ over all
$(i_{1},\ldots,i_{n}) \in Supp(f)$.
The convex hull of the set $Supp(f)$, denoted $P_{f}$, is known as the
{\em Newton polytope} of $f$.

The following lemma was observed by Ostrowski \cite{Ost21} in 1921 (see also
\cite[Theorem VI, p.\ 226]{Ost75}).

\begin{lemma}\label{folk}
Let $f,g,h \in K[X_{1},\ldots,X_{n}]$ with
$f=gh$. Then $P_{f} = P_{g}+P_{h}$.
\end{lemma}

An {\em integral polytope} is a polytope whose vertices have integer
coordinates, and we say that an integral polytope is {\em integrally
decomposable}, or simply {\em decomposable}, if it
can be written as a Minkowski sum of two integral polytopes,
each of which has more than one point. A summand in an integral
decomposition is called an {\em integral summand}.
We say an integral polytope is {\em integrally indecomposable}, or simply 
{\em indecomposable}, if it is not decomposable. 
The Newton polytope of a polynomial
is certainly integral and if the polynomial
factors into two polynomials each of which has at least two terms,
then by Lemma \ref{folk} its Newton polytope must be decomposable.
Thus we have the following simple irreducibility criterion from \cite{Gao99}.

\begin{corollary}[Irreducibility Criterion]\label{Gao}
Let $f\in K[X_{1},\ldots,X_{n}]$ with $f$ not divisible by
any $X_{i}$ for $1\leq i\leq n$. If the Newton polytope
of $f$ is integrally indecomposable, then $f$ is absolutely irreducible.
\end{corollary}

In Section \ref{sec-indecomp}, we shall discuss in more detail 
constructions of indecomposable polytopes and show 
how to get indecomposable polytopes
of high dimension from those of lower dimensions.  From these 
indecomposable polytopes one can easily give explicitly
many infinite families of polynomials which are absolutely
irreducible when considered over any field.

\subsection{Relevant computational problems}

 From a computational point of view, the following
problem is of interest.

\begin{problem}
Given an integral polytope, say as its list of vertices,
decide whether it is integrally indecomposable.
\end{problem}

This problem is not only pertinent to the study of polynomial
factorization, but is a natural problem to consider and as such may
be useful in other applications.
Here the input size is the length of the binary representation of the
coordinates of the vertices. 
Note that in  our applications the polytope will be
presented as the convex hull of a set of integral points. 
There is a large literature
on computing the convex hull of any finite set of points in $\R^n$;
see \cite[pages 361--375]{GO97}.
In particular, the convex hull of $t$ points in a plane can be
computed in time $O(t \log t)$ \cite{Gra72}.
Any of these algorithms can be used to
compute the vertices of the Newton polytope of a given polynomial and
we shall ignore this computational problem in the presentation of our
algorithms.
  
As mentioned before, the above problem is NP-complete, thus we shall be 
contented with algorithms that are ``efficient'' in terms
of some more generous measure, say the volume of polytopes.
In Section \ref{sec-decomp} we give such an algorithm for
polytopes in $\R^2$ and we also present a heuristic algorithm for
higher dimensional polytopes which uses randomization. 
It is an open problem to develop an ``efficient'' deterministic or even
randomized algorithm for testing general integral polytopes for
indecomposability.

For a decomposable integral polytope, it is desirable to find
all of its integral summands. Here we should identify polytopes
that are translations of each other. 

\begin{problem}
Given an integral polytope, say as its list of vertices,
find all of its integral summands.
\end{problem}

Again, this problem seems hard, but we shall give in Section \ref{sec-decomp}
an algorithm for polytopes of dimension two which is ``best possible''
in the sense that the running time is linearly related to the number of
decompositions.

\subsection{Some preliminary results}

We shall need more properties of the Minkowski sum.
The next result from \cite{Gao99} describes how the faces decompose 
in a Minkowski sum of polytopes;
for its proof, see Ewald~\cite[Theorem 1.5]{Ewa96},
Gr\"{u}nbaum~\cite[Theorem 1, p.\ 317]{Gru67}, 
or Schneider~\cite[Theorem 1.7.5]{Sch93}.

\begin{lemma}\label{lem3.2}
Let $P=Q+R$ where $Q$ and $R$ are polytopes in $\R^n$. Then
\begin{enumerate}
\item[(a)] Each face of $P$ is a Minkowski sum of unique faces of $Q$ and $R$.
\item[(b)] Let  $P_1$ be any face of $P$ and $c_1,\ldots, c_k$
all of its vertices. Suppose that $c_i=a_i + b_i$ where $a_i \in Q$
and $b_i \in R$ for $1 \leq i \leq k$. Let
$$Q_1 = \conv(a_1,\ldots, a_k), \ \ \
R_1 = \conv(b_1,\ldots, b_k). $$
Then $Q_1$ and $R_1$ are faces of $Q$ and $R$, respectively,
and $P_1 = Q_1 + R_1$.
\end{enumerate}
\end{lemma}

A polytope of dimension two is called a {\em polygon}.
(We refrain from using the term Newton polygon for a $2$-dimensional
Newton polytope as in number theory this term is used to refer to the
lower boundary of the ``Newton polyhedron'' of certain power series.)
The only proper faces of a polygon are its vertices and edges.
For polygons, the above lemma can be rephrased as follows.

\begin{corollary}\label{lemfaces}
Let $P,\,Q$ and $R$ be convex polygons (in $\R^n$) with $P=Q+R$.
Then every edge of $P$ decomposes uniquely as the sum of an edge of $Q$ and an
edge of $R$, possibly one of them being a point.
Conversely, any edge of $Q$ or $R$ is a summand of
exactly one edge of $P$.
\end{corollary}

\section{Indecomposable polytopes}\label{sec-indecomp}

First of all, we mention the following two constructions of 
indecomposable polytopes from \cite{Gao99}. 

\begin{theorem}\label{thm3.1}
Let $Q$ be any integral polytope in $\R^n$ contained
in a hyperplane $H$ and $v \in \R^n$ an integral point
lying outside of $H$.
Suppose that $v_1,\ldots, v_k$ are all the vertices of $Q$.
Then the polytope $\conv(v,Q)$ is integrally indecomposable iff
$$\gcd(v-v_1, \ldots, v-v_k)=1.$$
\end{theorem}
Here and hereafter the $\gcd$ of a collection of integral vectors 
is defined to be the $\gcd$ of all their coordinates together.

\begin{theorem}\label{thm8}
Let $Q$ be an indecomposable integral polytope in $\R^n$ that is contained
in a hyperplane $H$ and has at least two points, and let $v \in \R^n$ be
a point (not necessarily integral)
lying outside of $H$. Let $S$ be any set of integral points in the
polytope $\conv(v,Q)$. Then the polytope $\conv(S,Q)$ is integrally
indecomposable.
\end{theorem}

The first construction shows that an integral line
segment $\conv(v_0,v_1)$
is indecomposable iff $\gcd(v_0-v_1)=1$, and an integral triangle
$\conv(v_0,v_1,v_2)$ is integrally indecomposable iff
$\gcd(v_0-v_1,v_0-v_2)=1$. The second construction
gives many indecomposable polygons with more than three edges.
These two constructions can be used iteratively to get indecomposable
polytopes of any higher dimension. 

In the following, we give a new construction based on a projection.
Intuitively, one hopes that if a projection of a polytope is
indecomposable then the polytope is indecomposable itself. Unfortunately,
this is not true in general; consider for example a square and project it
along one of its edges. The following lemma, however, gives a sufficient
condition. We say that a linear map $\pi: \R^n \longrightarrow \R^m$
is {\em integral} if it maps integral points in $\R^n$ to integral
points in $\R^m$. It is straightforward to see that the image of
any integral polytope under an integral linear map is still an integral polytope.

\begin{lemma}\label{lem9}
Let $P$ be any integral polytope in $\R^n$ and
$\pi: \R^n \longrightarrow \R^m$ any integral linear map. If $\pi(P)$ 
is integrally indecomposable and each vertex of $\pi(P)$ has only one
preimage in $P$ then $P$ must be integrally indecomposable.
\end{lemma}
\begin{proof}
It suffices to show that $\pi(P)$ is decomposable if  $P$ is decomposable.
Suppose that $P= A+B$ for some integral polytopes $A$ and $B$
in $\R^n$ each with at least two points.  
Then $\pi(P) = \pi(A) + \pi(B)$.
We need to show that both $\pi(A)$ and $\pi(B)$ have at least two points.
Suppose otherwise, say $\pi(A)$ has only one point.
Let $w_0$ be any vertex of $P$ such that $\pi(w_0)$ is a vertex of $\pi(P)$.
Since $P=A+B$, there are unique vertices $u_0 \in A$ and $v_0 \in B$
such that $w_0 = u_0 + v_0$. As $A$ has at least two points, it has another
vertex $u_1$ such that $u_0u_1$ is one of its edges. Then,
by Lemma~\ref{lem3.2}, 
$P$ has an edge $w_0w_1$
that starts at $w_0$ and is parallel to $u_0u_1$ where $w_1$ is a
vertex of $P$ 
different from $w_0$.
 The latter property implies that $w_1 - w_0 = t(u_1-u_0)$
for some real number $t$. Hence
$$ \pi(w_1) - \pi(w_0) = \pi(w_1-w_0) = \pi(t(u_1 - u_0))
= t(\pi(u_1) - \pi(u_0)) = 0,$$
as $\pi(A)$ has only one point and $u_1,u_0 \in A$.  
This means that $\pi$ maps two vertices of $P$ to one vertex of $\pi(P)$,
contradicting our assumption.
\end{proof}

\begin{corollary}
Let $P$ be any integral polytope in $\R^n$ and 
$\pi: \R^n \longrightarrow \R^m$
any integral linear map that is injective on the vertices of $P$.
If $\pi(P)$ is integrally indecomposable then so must be $P$.
\end{corollary}

\begin{theorem}\label{thm11}
Let $Q$ be any integrally indecomposable polytope in $\R^m$ and
$\pi: \R^n \longrightarrow \R^m$ any integral linear map.
Let $S$ be any set of integral points in $\pi^{-1}(Q)$ having
exactly one point in $\pi^{-1}(v)$ for each vertex $v$ of $Q$.
Then the polytope $\conv(S)$ in $\R^n$ is integrally  indecomposable. 
\end{theorem}

\begin{proof} It follows directly from Lemma~\ref{lem9}.
\end{proof}

\noindent{\bf Remark}. Theorem~\ref{thm8} can be viewed as
a special case of Theorem~\ref{thm11} in the case that
$Q$ has sufficiently many
integral points besides its vertices, since it seems likely 
that there is an integral linear map that  projects integral
points in the cone $\conv(v,Q)$ to integral points in its base $Q$.
Such a projection is impossible if $Q$ has no integral points other
than its vertices.

In concluding this section, we would like to discuss
the relationship of integral decomposability with
a different concept of decomposability of polytopes defined  in 
Gr\"{u}nbaum~\cite[Chapter 15]{Gru67}.
Let $P,Q$ be polytopes in $\R^n$ (not necessarily integral).
We say that $Q$ is {\em homothetic} to $P$
if there is a real number $t \geq 0$ and a vector $a \in \R^n$ such that
$$ Q = t P + a = \{ t b + a: b \in P\}.$$
A polytope $P$ is called {\em homothetically indecomposable}
if it is the case that whenever
 $P=P_1 + P_2$ for any polytopes $P_1$ and $P_2$, then
$P_1$ or $P_2$ is homothetic to $P$. Otherwise, $P$ is called
{\em homothetically decomposable}.
Indecomposable polytopes in this sense have been extensively studied
in the literature \cite{Gal54,Kal82,McM87,Mey74,She63,Smi86,Smi87}.

Homothetic decomposability is not directly comparable with
integral decomposability.
On the one hand, the only 
homothetically indecomposable polytopes in the plane are 
line segments and triangles so any polygon with more than $3$ edges
is homothetically decomposable \cite{Gru67,Sch93}.
On the other hand, we saw above that some triangles can be
integrally decomposable and many polygons with more than $3$ edges
are integrally indecomposable!
The next result, however,  shows that homothetic indecomposability 
implies integral indecomposability under a simple condition.

\begin{proposition}\label{homo}
Let $Q$ be an integral polytope in $\R^n$ with vertices 
$v_i$, where $0 \leq i \leq k$. If $Q$ is homothetically 
indecomposable and
$$ \gcd(v_0-v_1, \cdots, v_0 - v_k) =1,$$
then $Q$ is integrally indecomposable.
\end{proposition}
\begin{proof}
Suppose that $Q = T + S$ for some integral polytopes $T$ and $S$.
Then $T$ or $S$ is homothetic to $Q$, say $T$.
This means that there is a real number $r\geq 0$ and $a\in \R^n$ such that
$T= r Q + a$.
Hence the vertices of $T$ are 
$$ u_i :=r v_i + a, \ \  i=0, 1, \ldots, k.$$
Since $T$ is integral, all the vertices $u_0, u_1, \ldots, u_k$
are integral and in particular
$$ u_0 - u_i = r (v_0-v_i), \ \ i=1, \ldots, k$$
are integral. So $r$ must be a rational number and the denominator
of $r$  divides $\gcd(v_0-v_1, \cdots, v_0 - v_k) =1$;
hence $r$ is an integer. As $0 \leq r \leq 1$, we have $r=0$ or $1$.
In either case, $T$ is a trivial summand of $Q$.
Therefore $Q$ is integrally indecomposable.
\end{proof}

By the above theorem, the homothetically indecomposable polytopes 
constructed in \cite{Gal54,Kal82,McM87,Mey74,Smi86,Smi87} give many
integrally indecomposable polytopes.

\section{Decomposing polytopes} \label{sec-decomp}

In this section we present our algorithms for both testing polytopes for
indecomposability and constructing summands of polytopes. We restrict
our attention to polygons in Section \ref{sec-polygon}
before considering the more general case in Section \ref{sec-polytope}.

\subsection{Polygons} \label{sec-polygon}

Given a convex
polygon in the Euclidean plane, one may form a finite
sequence of vectors associated with it as follows. Let
$v_{0},v_{1},\ldots,v_{m-1}$ be the vertices of the polygon
ordered cyclically in a clockwise direction.
The edges of $P$ are represented by the
vectors $E_i= v_{i}-v_{i-1} = (a_{i},b_{i})$ for $1 \leq i\leq m$, 
where $a_i, b_i \in \Z$ and the indices are taken modulo $m$. 
We call each $E_i$ an {\em edge vector}. 
A vector $v=(a,b) \in \Z^2$ is called a {\em primitive vector}
if $\gcd(a,b)=1$. 
Let $n_{i} = \gcd{(a_{i},b_{i})}$ and define $e_{i} = (a_i/n_i, b_i/n_i)$.
Then $E_i = n_i e_i$ where $e_i$ is a primitive vector, $1 \leq i\leq m$.
Each edge $E_{i}$ contains precisely $n_{i}+1$ integral points 
including its end points.  The sequence of vectors
$\{n_{i}e_{i}\}_{1\leq i\leq m}$, which we call the {\em edge sequence}
or a {\em polygonal sequence},
uniquely identifies the polygon up to translation determined by $v_0$, and will
be the input to our polygon decomposition algorithm.  It will be
convenient to identify sequences with those obtained by extending the
sequence by inserting an arbitrary number of zero vectors.  We may
thus assume that the edge sequence of a summand of a polygon $P$ has
the same length as that of $P$.  As the boundary of the polygon is a
closed path, we have that $\sum_{1\leq i\leq m} n_{i}e_{i}= (0,0)$.

\begin{lemma}\label{lemseq}
Let $P$ be a polygon with edge sequence $\{n_{i}e_{i}\}_{1\leq i \leq m}$
where $e_i\in \Z^2$ are primitive vectors.
Then an integral polygon is a summand of $P$ iff its edge sequence is of
the form $\{k_{i}e_{i}\}_{1\leq i \leq m}$, $0\leq k_i \leq n_i$, with
$\sum_{1\leq i\leq m} k_{i}e_i =(0,0)$.
\end{lemma}

\begin{proof} Let $\{e^{\prime}_{i}\}_{1\leq i\leq m}$ be the edge
sequence of an integral summand $Q$ of $P$.  
By the final statement in Corollary \ref{lemfaces},
each edge of $Q$ occurs as the summand of some edge $ne$ of $P$ where
$e$ is a primitive vector, and it is easily seen that its corresponding
edge vector must be of the form $ke$ with $0 \leq k \leq n$.
The sum is zero simply because the boundary of $Q$ is a closed path.
Conversely, any sequence of this form will determine a closed path. 
Since  $\{n_{i}e_{i}\}_{1\leq i \leq m}$
is a polygonal sequence, $\{k_{i}e_{i}\}_{1\leq i \leq m}$ must define
the boundary of a convex polygon.
It will be a summand of $P$, with the other summand 
having edge sequence $\{(n_{i}-k_{i})e_{i}\}_{1\leq i \leq m}$.
\end{proof}

Given as input a sequence of edge vectors $\{n_{i}e_{i}\}_{1\leq i\leq m}$
of a polygon $P$, our polygon
decomposition algorithm will check for the existence of a sequence
of integers  $k_{i}$ with $0\leq k_{i} \leq n_{i}$, $1\leq i \leq m$,  
such that $\sum_{1\leq i \leq m} k_{i}e_{i} = (0,0)$, $k_m \ne n_m$,
and not all $k_i=0$.
(If $P$ is decomposable then at least one of its summands has $k_m\ne n_m$.)
Thus the decision problem underlying our algorithm is

\vspace{0.3cm}

{\sc Polygon Decomposability (PolyDecomp)}\\
Input: The egde sequence $\{n_{i}e_{i}\}_{1\leq i\leq m}$ 
 of an integral convex polygon $P$.\\
Question: Does $P$ have a proper integral decomposition?\\

\vspace{0.3cm}

The input size of an instance of this problem is $O(m(\log N + \log E))$
where $N=\max\{n_1,\ldots, n_m\}$ and $E$ the maximum of absolute values
of the coordinates of $e_i$, $1 \leq i \leq m$. 
The next result puts the difficulty of this problem in context.

\begin{proposition}
{\sc PolyDecomp} is {\em NP}-complete.
\end{proposition}

\begin{proof}
Certainly the language associated with {\sc PolyDecomp} lies in NP as we
may use a proper decomposition of $P$ to verify membership of the language.
We give a polynomial reduction of {\sc Partition} to {\sc PolyDecomp}
which proves, since {\sc Partition} is NP-complete \cite{GJ79}, that {\sc
PolyDecomp} is NP-complete.

Recall that the input to {\sc Partition} is a sequence $\{s_{i}\}_{1
\leq i \leq m}$ of positive integers which we may take to be non-decreasing.
Thus $s_{1} \leq s_{2} \leq \ldots \leq s_{m}$. Let $t = \sum_{1 \leq i
\leq m} s_{i}$. The question in {\sc Partition} is whether there is a
subsequence of $\{s_{i}\}$ with sum $t/2$. Observe that we may assume
that $t$ is even, for otherwise the question is easily
answered. Consider now the following instance of {\sc PolyDecomp}:
the edge sequence 
\[(s_{1},1),(s_{2},1),\ldots
(s_{m},1), m(0,-1), (-t/2,-1), (-t/2,1)\]
where all $n_i=1$. Firstly,
it is easy to check that this is indeed a polygonal sequence. Secondly, 
any polygon associated with the polygonal sequence has a proper decomposition if
and only if the sequence $\{s_{i}\}_{1 \leq i \leq m}$ has a
subsequence with sum $t/2$. Thus we have a polynomial reduction,
which completes the proof.
\end{proof}

Since it is widely believed that NP $\ne$ P, it seems unreasonable to attempt
to find a genuinely efficient algorithm for solving {\sc PolyDecomp};
however, we shall present an algorithm below whose running time is
polynomial in the length of the sides of the polygon rather than 
the logarithm of the lengths. In the parlance of \cite{GJ79}, 
this is an example of a ``pseudopolynomial-time'' algorithm. 
In Section \ref{sec-appl} we shall indicate how this algorithm may be used
to test bivariate polynomials for
absolute irreducibility; the algorithm 
thus obtained is efficient in terms of the total degree of the
polynomial, rather than the number of non-zero terms. Thus the distinction
between genuinely efficient algorithms for deciding polytope
decomposability and ``pseudopolynomial-time'' algorithms is mirrored
to a certain extent in that between efficient algorithms
for polynomials in terms of their sparse and dense representations.

\begin{algorithm}\label{alg-polygon}
{\bf (PolyDecomp)}\\
{\it Input:} The edge sequence $\{n_{i}e_{i}\}_{1 \leq i\leq m}$
of an integral convex polygon $P$ starting at a vertex $v_{0}$ 
where $e_i \in \Z^2$ are primitive vectors.\\
{\it Output:} Whether $P$ is decomposable.
\vspace{0.3cm}

\hspace{-\parindent}{\it Step 1:}
Compute the set {\itshape IP} of all the integral points in $P$, 
and set $A_0=\emptyset$.
\vspace{0.3cm}

\hspace{-\parindent}{\it Step 2:} 
For $i$ from $1$ up to $m-1$, compute the set $A_i$ of points in
{\itshape IP} that are
reachable via the vectors $e_1, \ldots, e_i$:
\begin{enumerate}
\item[2.1] For each $0 < k \leq n_i$, if $v_0+ k e_i \in \mbox{\itshape
IP}$ then  add it to $A_i$;
\item[2.2] For each $u \in A_{i-1}$ and $0 \leq k \leq n_i$, if $u+k
e_i \in \mbox{\itshape IP}$
      then  add it to $A_i$.
\end{enumerate}

\hspace{-\parindent}{\it Step 3:} Compute the last set $A_m$:
 For each $u \in A_{m-1}$ and $0 \leq k < n_m$, if $u+k e_m \in \mbox{
\itshape IP}$ then  add it to $A_m$.
\vspace{0.3cm}

\hspace{-\parindent}{\it Step 4:} Return ``Indecomposable'' if $v_0 \not \in A_m$ and ``Decomposable'' otherwise.
\end{algorithm}

\begin{theorem}
The above algorithm decides decomposability correctly in $O(tmN)$ vector 
operations where $t$ is the number of integral points in $P$, $m$ the number
of edges and $N$ the maximum number of integral points on an edge.
\end{theorem}

\begin{proof}
(Note that by a vector operation we mean adding two vectors, multiplying a vector
by a scalar, or adjoining a point to a set.)
The running time is easy to see as each set $A_i$ has size at most
$t$.
Also, the set {\itshape IP} can be computed in time $O(t)$: since the
edge sequence is already given one can enumerate points in $P$ by
scanning vertical line segments starting at $v_0$. One need only keep
track of the top and bottom edges when moving the lines
(i.e. increasing $x$ values) and the edges tell us the range of the
$y$ value for any given $x$ value. (Note that $t$ itself can
be bound in terms of $m$, $N$ and the area of the polygon using
Pick's formula \cite[page 139]{GO97}.)

To prove the correctness, observe that
all the points in $A_m$ are of the form
$v_0 + \sum_{i=1}^{m} k_i e_i$, $0 \leq k_i \leq n_i$.
Step 2.1 ensures that $k_i \ne 0$ for some $i<m$  and
Step 3 insists  that $k_m < n_m$
(note that $v_0 + k e_m \not\in \mbox{\itshape IP}$ for all $k>0$).
If one of the points in $A_m$ is equal to $v_0$ then
$\sum_{i=1}^{m} k_i e_i = (0,0)$, and so the sequence
$\{k_ie_i\}$ forms the edge sequence of a proper integral summand of $P$.
On the other hand, for any proper integral summand $Q$ of $P$,
$Q$ can be ``slid'' into $P$ at $v_0$, that is, $Q$ can be translated
so that $v_0$ is a vertex of $Q$ and $Q$ lies inside $P$. Hence
all the vertices of $Q$ must lie in $P$ and thus in {\itshape IP}. Consequently
its edge sequence will be detected by our algorithm.
\end{proof}

We next give a simple generalisation of the above algorithm which not only 
outputs the number of proper decompositions of the polygon, but also
outputs an array. The array may then be used to recover all decompositions,
a single ``recovery'' requiring linear time. Thus the total time taken
to recover all decompositions is essentially linearly related to the
number of decompositions. This is the best that one can expect;
however, it does not yield a ``pseudopolynomial-time'' algorithm as
the number of decompositions may be exponential in the area of the polygon. 
For example, consider the polygon with edge sequence 
$$(1,1), (2,1), \ldots , (m,1), m(0,-1), t(-1,0)$$
where $t=(m+1)m/2$.  The polygon has area less than
$1^2 + 2^2 + \cdots + m^2 = O(m^3)$ while the number of integral
summands is exactly $2^m$.

\begin{algorithm}\label{alg-polynum}
{\bf (PolyDecompNum)}\\
{\it Input:} The edge sequence $\{n_{i}e_{i}\}_{1 \leq i\leq m}$
of an integral convex polygon $P$ starting at a vertex $v_{0}$ 
where $e_i\in \Z^2$ are primitive vectors.\\
{\it Output:} The number of integral summands of $P$ including
the trivial ones, and an array $A$. Each cell in $A$ contains
a pair $(u,S)$ where
$u$ is a non-negative integer and $S$ is a subset of $\{(k,i): 1 \leq 
k \leq n_{i}, 1 \leq i \leq m\}$.
\vspace{0.3cm}

\hspace{-\parindent}{\it Step 1:}
Compute the set {\itshape IP} of all the integral points in $P$ 
(so $v_0 \in \mbox{\itshape IP}$); say {\itshape IP}
has $t$ points.
Initialize a $t$-array $A_0$ indexed by the points in {\itshape IP}.
Set $A_0[v] := 
(0,\emptyset)$ for all $v \in \mbox{\itshape IP}$
except the cell $A_{0}[v_0]$ which is set to $(1,\emptyset)$.
\vspace{0.3cm}

\hspace{-\parindent}{\it Step 2:} 
For $i$ from $1$ up to $m$, compute the $t$-array $A_i$ from $A_{i-1}$:
\begin{enumerate}
\item[2.1] First copy the contents of all the cells of $A_{i-1}$ into $A_i$
  (this step is for $k=0$).
\item[2.2] For each $v \in \mbox{\itshape IP}$ with the first number of the cell $A_{i-1}[v]$ nonzero,
and for each $0 < k \leq n_i$, if $v^{\prime} = v+k e_i \in
\mbox{\itshape IP}$ then update the cell $A_i[v^{\prime}]$ as follows: if $(u_1, S_1)$ is the value
of $A_{i-1}[v]$ and $(u_2, S_2)$ the current value of $A_i[v^{\prime}]$ then the new
value of $A_i[v^{\prime}]$ is $(u_1+u_2, S_2 \cup \{(k,i)\})$.  
\end{enumerate}
\vspace{0.3cm}

\hspace{-\parindent}{\it Step 3:} Return the number $u$ and the array $A=A_m$,
where $(u,S)$ is the content of cell $A_{m}[v_0]$.
\end{algorithm}

\begin{theorem}
The integer output by Algorithm \ref{alg-polynum}
is the total number of
integral summands of the polygon $P$.
\end{theorem}

\begin{proof}
Supposing $v=v_0+ k_1e_1 + \cdots + k_ie_i$, we may view the vector sum 
as a path from $v_0$ to $v$, so the number of such paths is equal to 
the sum of the numbers of paths from $v_0$ to $v-ke_i$ for $0\leq k\leq n_i$, 
using $e_1, \ldots, e_{i-1}$. Hence the numbers of paths can be computed
iteratively as described in the algorithm:
the number $u$ in $A_i[v]$ records the number of paths from $v_0$ to
$v$ using $e_1, \ldots, e_{i}$ and the set $S$ records all the pairs $(k,j)$,
$j \leq i$, for which a path reaches $v$ with its last edge being $ke_j$ 
with $k>0$. 
Thus the integer in cell $A_{m}[v_{0}]$ is the total number of closed
paths $\sum_{1 \leq i \leq m} k_{i}e_{i}$ starting at $v_{0}$.
By Lemma \ref{lemseq} this is the number of integral summands of $P$.
\end{proof}

The significance of the array $A$ output by the algorithm is that it
may be used to recover all decompositions of the polygon $P$. We show
how a single decomposition can be recovered: Suppose the cell
$A[v_0]$ contains the pair $(u,S)$. Choose any $(k,i) \in S$. The
line segment $ke_i$  will be the ``final edge''
(counting clockwise) in our
summand of $P$. Let $(u^{\prime},S^{\prime})$ be the contents of cell
$B[v_0 - ke_i]$. Pick any $(k^{\prime},i^{\prime}) \in S^{\prime}$ with
 $i^{\prime} < i$. The line segment 
$k^{\prime}e_{i^{\prime}}$ will be the ``penultimate edge'' in
our summand of $P$. We continue in this way, and as our sequence of
$i$'s is decreasing we shall eventually return to the cell
$A[v_0]$. At that point we will have recovered one summand in a decomposition
of $P$.

With regard to the running time, each cell in the array can be updated
at most $mN$ times, thus the running time is $O(tmN)$ ``cell
updates''. The data in each cell is a pair $(u,S)$ where $S$ is a set
of size at most $mN$ and $u$ an integer less than $N^{m}$ (an upper
bound on the number of summands). Updating the integer $u$
involves integer addition and this has a bit complexity of
$O(\log{N^{m}}) = O(m\log N)$. Updating the set $S$ simply involves unioning
it with an element $(k,i)$. Ignoring logarithmic factors, we can
consider this a single bit operation. Thus the running time 
of {\bf PolyDecompNum} is 
$O(tm^{2}N)$ bit operations, ignoring logarithmic factors.

\subsection{Higher dimensional polytopes} \label{sec-polytope}

The problem of testing higher dimensional polytopes for
decomposability appears to be significantly more difficult.
Certainly it is NP-complete as it includes that of polygons as a special
case.  It would be interesting to investigate whether
this problem was ``strongly NP-complete'' in the sense of \cite{GJ79};
this essentially means that the problem remains ``NP-complete'' when one
bounds running time by the lengths, instead of logarithm of the lengths,
of the edge vectors.
If this more general problem is ``strongly NP-complete'' then it is
unlikely there is an algorithm for determining whether a convex
polytope of arbitrary dimension is indecomposable whose running time
is polynomial in terms of the volume of the polytope. 

In this section, we present a heuristic ``randomized algorithm''
based on the projections considered in 
Lemma~\ref{lem9}. The algorithm has running time polynomial in the
lengths of the edges of the
polytope, thus is ``efficient'' in the sense which we
have been considering. 
The idea is to choose a random integral linear map that
projects a polytope into a polygon in a plane and then test the decomposability
of the polygon. If the polygon is indecomposable and the condition of 
Lemma~\ref{lem9} is satisfied then the original polytope
is indecomposable. We will show that the condition of  Lemma~\ref{lem9}
is always satisfied with high probability, but we do not know how to prove
a good bound on the probability that the projected polygon be indecomposable
when the original polytope is indecomposable.

We now describe the details of our algorithm.
Let $S \subset \R^n$ be any finite set of integral points, 
which will be the input to our algorithm, and  $P=\conv(S)$.
We want to decide whether $P$ is integrally indecomposable.
Note that $P$ can be computed from $S$ by any of the
algorithms in \cite{GO97,Gra72}; however, 
our algorithm does not require that the vertices, which are all in $S$,
of $P$ be known in advance but detects them automatically. 
This is because the points of $S$ 
that are mapped to vertices of a polygon will be vertices of $P$,
provided each vertex of the polygon has only one preimage in $S$.

To describe a projection, we write points in $\R^n$ as column
vectors, so a set $S$ of $\ell$ points can be represented as an $n\times \ell$
matrix where each column stands for a point; for convenience, we still
denote the matrix by $S$.  As the points in $S$ are distinct so are
the columns of $S$.  Let $u,v \in \R^n$ be two integral points.
Then for any point $w\in \R^n$, the matrix-vector product $(u,v)^t w$
can be viewed as a point in $\R^2$. This defines an integral projection
$\pi$ from $\R^n$ into $\R^2$ and  
\begin{equation}\label{pp}
 (u,v)^t S
\end{equation}
is the  image of $S$ under $\pi$ in $\R^2$.
The polygon defined by the convex hull of the points in (\ref{pp}) is called
the {\em shadow polygon}, or simply {\em shadow},
of $P$ projected by $u$ and $v$.
The next lemma from \cite{Gao99a} arises in a different context and
tells us how likely it is that the projection is injective on the set $S$;
its proof is straightforward.

\begin{lemma}\label{lem2.5}
Let $S$ be an $n\times \ell$ matrix over a field
with no repeated columns and let $K$ be any subset of
cardinality $k$ of the same field. Pick $u_i \in K$ randomly and
independently, $1 \leq i \leq n$, and let
$$(a_1, \cdots, a_\ell) = (u_1,\cdots, u_n)S.$$
Then with  probability at least $1-\frac{\ell (\ell-1)}{2k}$ the entries
$a_1, \ldots, a_\ell$ are distinct.
\end{lemma}
 
Now let $K=\{-\ell^2, \ldots, -1, 0 , 1, \ldots, \ell^2\}$ which has 
$k=2\ell^2+1$ integers. If we choose the entries of 
$u$ and $v$ from $K$ at random and independently, then 
with probability at least $3/4$ the points in (\ref{pp}) are distinct,
so the condition in Lemma \ref{lem9} is satisfied, i.e., each vertex
of the shadow has only one preimage in $S$.
This probability can be increased arbitrarily close to $1$ if one increases
the size of the set $K$.

\begin{algorithm}\label{alg-polytope}
{\bf (PolytopeDecomp)}\\
{\it Input:} A finite set $S$ of integral points in $\R^n$.
\vspace{0.3cm}

\hspace{-\parindent}{\it Output:} 
``Indecomposable''or ``Failure''; the first case means
that the polytope $P=\conv(S)$ is proved to be indecomposable 
while the latter means the decomposability of $P$ is not decided.   
\vspace{0.3cm}

\hspace{-\parindent}{\it Step 0:} Form the points in $S$ as an $n\times \ell$
matrix, still denoted by $S$, where $\ell$ is the cardinality of $S$ and 
each column represents  a point.  Fix a set $K$ of small integers.
\vspace{0.3cm}

\hspace{-\parindent}{\it Step 1:}
Pick two vectors $u,v \in K^n$ randomly and compute
the projection $(u,v)^t S =(a_1, \ldots, a_\ell)$ where $a_i \in \Z^2$.
\vspace{0.3cm}

\hspace{-\parindent}{\it Step 2:}
Compute the vertices, say $v_1, \ldots, v_m$ in a clockwise direction,
 of the convex polygon 
defined by the points $a_1, \ldots, a_\ell$.  If more than two points
of $S$ are mapped to  one of the vertices $v_i$'s, then output ``Failure''
and stop here.
\vspace{0.3cm}

\hspace{-\parindent}{\it Step 3:}
Compute $E_i = v_i - v_{i-1} = n_i e_i$ where $n_i$ is a positive integer
and $e_i$ is a primitive vector, $1 \leq i \leq m$.

\vspace{0.3cm}
\hspace{-\parindent}{\it Step 4:}
Input the edge sequence $\{ n_i e_i\}$ to Algorithm 
{\bf PolyDecomp}. If the latter says ``Indecomposable''
then output ``Indecomposable'', otherwise output ``Failure''. 
\end{algorithm}

The correctness of this algorithm follows from our discussion
above. If $P$ is integrally decomposable then the algorithm will
always output ``Failure''.  It remains an open problem to determine
how likely it is that the algorithm will output ``Indecomposable'' if $P$ is 
integrally indecomposable.
It is possible that there are indecomposable polytopes whose
shadow polygons are always decomposable; for such
polytopes our algorithm will not work. We would be very interested in
seeing such examples.

On the other hand, it can be proved that most polytopes in 
$\R^n$, $n \geq 3$, are homothetically indecomposable 
\cite[Theorem 3.2.14, p152]{Sch93}. By Proposition \ref{homo},
we may expect that most integral polytopes are integrally indecomposable
so our algorithm  may detect most of them quickly.
It would be interesting to know how likely it is that
a random shadow polygon of
a random integral polytope (under some  probability distribution)
is indecomposable.

\section{Applications to polynomials} \label{sec-appl}

A direct application of Algorithm \ref{alg-polygon} in the light of
Corollary \ref{Gao} gives an algorithm for testing absolute
irreducibility of bivariate polynomials. One simply first checks
whether the input polynomial has any factors of the form $X_{i}$ and
if not computes the edge sequence of its Newton polytope, which can be done
in $O(t\log t)$ operations where $t$ is the number of nonzero terms
in the polynomial. 
Algorithm \ref{alg-polygon} may then be used to determine whether this polygon
is decomposable; if it is indecomposable then the polynomial must be
absolutely irreducible. In the case that the polygon is decomposable
the test is inconclusive. The running time of this algorithm is easily
checked to be $O(n^{3})$ where $n$ is the total degree of the
polynomial. A similar test based on Algorithm \ref{alg-polytope} may
be devised to test general multivariate polynomials for absolute
irreducibility where $S$ is taken to be the set of support vectors of the
polynomial to be tested. 

Certainly, this polytope approach cannot decide irreducibility of some
polynomials since it uses only their ``shapes'', i.e.  Newton
polytopes, and the coefficients do not come into play.  However, our
algorithm is extremely fast compared to the infallible algorithms in
\cite{Duv91,GK85,GC84,Kal85,Len85,Len87}, thus it may be used as a
pretest before applying the more expensive methods.  For random sparse
polynomials, their Newton polytopes may be viewed as random integral
polytopes. As we mentioned at the end of the last section, most
integral polytopes are expected to be indecomposable.  Hence the
``shapes'' of most polynomials are indecomposable, so our algorithm
can detect them quickly in most of the cases.  This means that our
polytope method should be particularly effective for random sparse
polynomials. The reader is referred to \cite{GL2000} for an
implementation of this algorithm which gives more precise details on
the range of applicability and effectiveness of the approach.

We finish by returning to the problem of Sturmfels quoted in 
Section \ref{sec-intro}. In this problem, one is given the list
of support vectors of a polynomial $f$ but the coefficients 
of $f$ are not specified. From the  support vectors, one can
compute their convex hull. So one is essentially given the
Newton polytope $P_f$ of $f$ with the requirement that the terms of
$f$ corresponding to the integral points of $P_f$ not on the given list
of support vectors must have zero coefficient. The question is how
such a polynomial $f$ factors in general? What are the Newton polytopes
and support vectors for the factors?

A natural approach is to find the set of all integral summands of
$P_f$, as this set
contains the Newton polytopes of all possible factors.  Each
summand may correspond to a factor of $f$, and if this is the case
then the set of integral points in the summand contains the
support vectors of the corresponding factor.  For bivariate
polynomials, one may find all integral summands by applying Algorithm
\ref{alg-polynum} and the method suggested immediately after it.  It
seems that most integral polytopes do not have many integral summands,
so our method is expected to be effective for random sparse polynomials. 
We would like to add that this method can be refined by taking into
account the possible factorizations of the univariate polynomials defined
by the edges of the polygon; however, we do not pursue this at present.

We should point out that some integral summands may not correspond
to any factor of $f$. For example, let 
$$f= (a + bX^n) + Y^m (c+dX^n) \in K[X,Y].$$
Its Newton polytope is a rectangle defined
by the support vectors $(0,0)$, $(0,n)$, $(m,0)$ and $(n,m)$.
This rectangle has $(n+1)(m+1)$ integral summands. But $f$ is
almost always absolutely irreducible except for a few cases!
(Absolute irreducibility of such $f$ under mild conditions may
be proved by substituting $X$ for $X + \alpha$, where $\alpha$
is a root of $a + bX^n$, and considering Newton polytopes once again.
Of course, $f$ may have a univariate factor, and will be reducible
if the field characteristic divides both $m$ and $n$, but these
are the only exceptions.)
Moreover, in general even when we find a summand $P_{g}$ of $P_{f}$ which
corresponds to a factor $g$ of the polynomial $f$ under consideration,
it may be the case that not all integral points in $P_{g}$ are
support vectors of $g$. We only know for sure that the vertices of
$P_{g}$ are among the support vectors of $g$.

Finally, we mention that deciding reducibility of
sparse polynomials can be considered a special case of
the above problem. Even though we have shown that
deciding decomposability of integral polytopes is NP-complete,
we still do not know whether deciding reducibility
is also NP-complete. The latter problem is not
even known for sparse univariate polynomials over finite fields.

\section{Conclusion}

The Newton polytope of a polynomial carries a lot of information
about its factors, and so it is fruitful to study algorithms
for deciding decomposability of integral polytopes and
for finding all the integral summands when they are decomposable.
For polygons, we showed that deciding decomposability is
NP-complete but gave a pseudo-polynomial time algorithm for
testing decomposability  and for constructing all possible decompositions.
For polytopes of dimension larger than two, we presented 
an indecomposability lemma based on projections, and this
lemma gives a heuristic method for testing
decomposability of polytopes in any dimension.
However, a rigorous analysis of this algorithm is still lacking.
It is also desirable to have an algorithm for finding
all the integral summands for  polytopes in arbitrary
dimensions. 
The corresponding problems for (sparse) polynomials are also
open: it is not even known whether deciding reducibility of
sparse polynomials is NP-complete.

\end{document}